\documentclass{gtart}
\usepackage{amssymb,amsmath,verbatim}

\input gtoutput
\volumenumber{2}\papernumber{5}\volumeyear{1998}
\pagenumbers{79}{101}\published{3 June 1998}
\shorttitle{Flag Manifolds and the Landweber-Novikov Algebra}  
\proposed{Haynes Miller}\seconded{Gunnar Carlsson, Ralph Cohen}
\received{23 October 1997}\revised{6 January 1998}
\accepted{1 June 1998}
\def\S{section }

\newtheorem{prop}[equation]{Proposition}

\newtheorem{thm}[equation]{Theorem}
\newtheorem{cor}[equation]{Corollary}
\newtheorem{lem}[equation]{Lemma}

\theoremstyle{definition}

\numberwithin{equation}{section}

\newcommand{\spandsp}{\mbox{$\qquad\text{and}\qquad$}}
\newcommand{\sands}{\mbox{$\quad\text{and}\quad$}}
\newcommand{\sptsp}[1]{\mbox{$\qquad\text{#1}\qquad$}}
\newcommand{\sts}[1]{\mbox{$\quad\text{#1}\quad$}}

\newcommand{\Hom}{\operatorname{Hom}}

\newcommand{\dtimes}{\mathop{\mbox{\Large$\times$}}}

\newcommand{\ahss}{Atiyah--Hirzebruch spectral sequence}

\newcommand{\CP}{\text{\it CP\/}}
\newcommand{\CPI}{\mbox{$\text{\it CP\/}^\infty$}}

\newcommand{\pd}{\text{\it Pd\/}}
\newcommand{\hd}{\text{\it Hd\/}}

\newcommand{\MU}{\text{\it MU\/}}

\newcommand{\DU}{\text{\it DU\/}}

\newcommand{\osdu}{\mbox{$\varOmega_*^{DU}$}}
\newcommand{\odus}{\mbox{$\varOmega^*_{DU}$}}
\newcommand{\osu}{\mbox{$\varOmega_*^{U}$}}
\newcommand{\ous}{\mbox{$\varOmega^*_{U}$}}
\newcommand{\MUMU}{\mbox{${\mathit{MU}}\wedge{\mathit{MU}}$}}

\newcommand{\BU}{\text{\it BU\/}}

\newcommand{\bZ}{\mathbb{Z}}
\newcommand{\bR}{\mathbb{R}}
\newcommand{\bC}{\mathbb{C}}

\newcommand{\cD}{\mathcal{D}}

\begin{document}
\bibliographystyle{plain}
\title{Flag Manifolds and the Landweber--Novikov Algebra}
\author{Victor M Buchstaber\\Nigel Ray}
\address{Department of Mathematics and Mechanics, Moscow State
University\\119899 Moscow, Russia}
\secondaddress{Department of Mathematics, University of Manchester\\
Manchester M13 9PL, England}
\asciiaddress{Department of Mathematics and Mechanics, Moscow State
University\\119899 Moscow, Russia\\
Department of Mathematics, University of Manchester\\
Manchester M13 9PL, England}

\email{buchstab@mech.math.msu.su\\nige@ma.man.ac.uk}



\begin{abstract} 

We investigate geometrical interpretations of various structure maps
associated with the Landweber--Novikov algebra $S^*$ and its integral
dual $S_*$. In particular, we study the coproduct and antipode in
$S_*$, together with the left and right actions of $S^*$ on $S_*$ which 
underly the construction of the quantum (or Drinfeld) double
$\mathcal{D}(S^*)$. We set our realizations in the context of double
complex cobordism, utilizing certain manifolds of bounded flags which
generalize complex projective space and may be canonically expressed
as toric varieties. We discuss their cell structure by analogy with
the classical Schubert decomposition, and detail the implications for
Poincar\'e duality with respect to double cobordism theory; these lead
directly to our main results for the Landweber--Novikov algebra.  

\end{abstract}

\asciiabstract{We investigate geometrical interpretations of various
structure maps associated with the Landweber-Novikov algebra S^* and
its integral dual S_*. In particular, we study the coproduct and
antipode in S_*, together with the left and right actions of S^*
on S_* which underly the construction of the quantum (or Drinfeld)
double D(S^*). We set our realizations in the context of
double complex cobordism, utilizing certain manifolds of bounded flags
which generalize complex projective space and may be canonically
expressed as toric varieties. We discuss their cell structure by
analogy with the classical Schubert decomposition, and detail the
implications for Poincare duality with respect to double cobordism
theory; these lead directly to our main results for the
Landweber-Novikov algebra.}

\maketitle

\section{Introduction}\label{intro} 

The Landweber--Novikov algebra $S^*$ was introduced in the 1960s as an
algebra of cohomology operations in complex cobordism theory, and was
subsequently described by Buchstaber and Shokurov \cite{bush:lna} in
terms of differential operators on a certain algebraic group. More
recently, both $S^*$ and its integral dual $S_*$ have been studied
from alternative viewpoints \cite{kashue:ccr}, \cite{no:vdh},
\cite{rasc:cmc}, reflecting the growth in popularity of Hopf algebras
throughout mathematics. Nevertheless, the interpretations have
remained predominately algebraic, although the underlying motivations
have ranged from theoretical physics to combinatorics.

Our purpose here is to provide a purely geometric description of
$S_*$, incorporating its structure maps and certain left and right
actions by $S^*$; the importance of the latter is their contribution
to the adjoint action, which figures prominently in the construction
of the quantum (or Drinfeld) double $\cD(S^*)$. We work in the context
of double complex cobordism, whose properties we have developed in
a preliminary article \cite{bura:dcq}. So far as we are
aware, double cobordism theories first appeared in \cite{ra:phd}, and
in the associated work \cite{rasw:oss}. To emphasize our geometric
intent we return to the notation of the 60s, and write bordism and
cobordism functors as $\varOmega_*(\;)$ and $\varOmega^*(\;)$
throughout. 

The realizations we seek are provided by a family of bounded flag
manifolds with various double $U$--structures. These manifolds were
originally constructed by Bott and Samelson \cite{bosa:atm} without
reference to flags or $U$--structures, and were introduced into complex
cobordism theory in \cite{ra:ocb}. We consider their algebraic
topology in detail, describing computations in bordism and cobordism
theory which provide the essential link with the Landweber--Novikov
algebra, and are related to the generalized Schubert calculus of
Bressler and Evens \cite{brev:scc}. These results appear to be of
independent interest and extend to the topological study of other
toric manifolds \cite{babe:ktt}, \cite{daja:cpc}, as well as being
related to Magyar's program \cite{ma:bsv} for the description of
arbitrary Bott--Samelson varieties in combinatorial terms. We hope
to record such extensions in a future work.

For readers who seek background information in algebra, combinatorics,
and geometry, we suggest the classic books by Kassel \cite{ka:qg},
Aigner \cite{ai:ct}, and Griffiths and Harris \cite{grha:pag}
respectively. 

We begin in \S\ref{dococo} by summarizing prerequisites and notation
connected with double complex cobordism, recalling the coefficient
ring $\osdu$ and its subalgebra $G_*$, together with the canonical
isomorphism which identifies them with the Hopf algebroid $A_*^U$ and
its sub-Hopf algebra $S_*$ respectively. In \S\ref{boflma} we study
the geometry and topology of the bounded flag manifolds $B(Z_{n+1})$,
describing their toric structure and introducing the posets of
subvarieties $X_Q$ which serve to desingularize their cells. In
\S\ref{nostdu} we define the basic $U$-- and double $U$--structures on
$X_Q$ which underlie the geometrical realization of $G_*$, and use
them to compute $\ous(X_Q)$ and $\osu(X_Q)$; the methods extend to
double cobordism, although several aspects of duality demand extra
care. We apply this material in \S\ref{ap} to calculate $\osdu$--theory
characteristic numbers of the $X_Q$, interpreting the results by means
of the calculus of \S\ref{boflma}. Under the canonical isomorphism,
realizations of the relevant structure maps for $S^*$ and $S_*$ follow
immediately.

We use the following notation and conventions without further comment.

Given a complex $m$--plane bundle $\xi$ over a finite CW complex, we
let $\xi^\perp$ denote the complementary $(n-m)$--plane bundle in some
suitably high-dimensional trivial bundle $\bC^n$.

We write $A_U^*$ for the algebra of complex cobordism operations, and
$A^U_*$ for its continuous dual
$\Hom_{\varOmega^U_*}(A_U^*,\varOmega^U_*)$, forcing us in turn to
write $S^*$ for the graded Landweber--Novikov algebra, and $S_*$ for
its dual $\Hom_\bZ(S^*,\bZ)$; neither of these notations is entirely
standard.

Several of our algebras are polynomial in variables such as $b_k$ of
grading $2k$, where $b_0$ is the identity. An additive basis is
therefore given by monomials of the form
$b_1^{\omega_1}b_2^{\omega_2}\dots b_n^{\omega_n}$, which we denote by
$b^\omega$, where $\omega$ is the sequence of nonnegative, eventually
zero integers $(\omega_1,\omega_2, \dots,\omega_n,0\dots)$. The set of
all such sequences forms an additive semigroup, and $b^\psi
b^\omega=b^{\psi+\omega}$. Given any $\omega$, we write $|\omega|$ for
$2\sum i\omega_i$, which is the grading of $b^\omega$.  We distinguish
the sequences $\epsilon(m)$, which have a single nonzero element 1 and
are defined by $b^{\epsilon(m)}=b_m$ for each integer $m\geq 1$. It is
often convenient to abbreviate the formal sum $\sum_{k\geq 0}b_k$ to
$b$, and write $(b)^n_k$ for the component of the $n$th power of $b$
in grading $2k$; negative values of $n$ are permissible.

When dualizing, we choose dual basis elements of the form $c_\omega$,
defined by $\langle c_\omega,b^\psi\rangle=\delta_{\omega,\psi}$; this
notation is designed to be consistent with our convention on gradings,
and to emphasize that the elements $c_\omega$ are not necessarily
monomials themselves.

The authors are indebted to many colleagues for enjoyable and
stimulating discussions which have contributed to this work. These
include Andrew Baker, Sara Billey, Francis Clarke, Fred Cohen, Sergei
Fomin, Sergei Gelfand, Christian Lenart, Peter Magyar, Haynes Miller,
Jack Morava, Sergei Novikov, and Neil Strickland.

\section{Double complex cobordism}\label{dococo}
 
In this section we summarize the appropriate parts of \cite{bura:dcq},
concerning the notation and conventions of double complex cobordism,
and operations and cooperations in the corresponding homology and
cohomology theories.  

Double complex cobordism is based on manifolds $M$ whose stable normal
bundles are equipped with a splitting of the form
$\nu\cong\nu_\ell\oplus\nu_r$. We refer to an equivalence class of such
splittings as a {\it double $U$--structure\/} $(\nu_\ell,\nu_r)$ on $M$,
writing $(M;\nu_\ell,\nu_r)$ if the manifold requires emphasis. It is
helpful to think of $\nu_\ell$ and $\nu_r$ as the {\it left\/} and {\it
right normal bundles\/} of the structure, respectively. We may follow
Stong \cite{st:nct} and Quillen \cite{qu:eps} in setting up the
corresponding bordism and cobordism functors geometrically, taking
necessary care with the double indexing inherent in the
splitting. Cartesian product ensures that $\odus(X)$ is a graded ring
for any space or spectrum $X$. Both functors admit an involution $\chi$
induced by interchanging the order of $\nu_\ell$ and $\nu_r$, and we
find it convenient to write $\chi(M)$ for $(M;\nu_r,\nu_\ell)$. The
coefficient ring $\osdu$ is the {\it double complex cobordism ring}.

We may recombine the left and right normal bundles to obtain a
forgetful homomorphism $\pi\colon\osdu(X)\rightarrow\osu(X)$;
conversely, we may interpret any standard $U$--structure as
either of the double $U$--structures $(\nu,0)$ or $(0,\nu)$,
thereby inducing multiplicative natural transformations $\iota_\ell$
and $\iota_r\colon\osu(X)\rightarrow\osdu(X)$, which are interchanged
by $\chi$. All these transformations have cohomological counterparts,
and the compositions $\pi\circ\iota_\ell$ and $\pi\circ\iota_r$
reduce to the identity. Given an element $\theta$ of $\osu(X)$ or
$\ous(X)$, we write $\iota_\ell(\theta)$ and $\iota_r(\theta)$ as 
$\theta_\ell$ and $\theta_r$ respectively.

From the homotopy theoretic viewpoint, it is convenient to work in any
of the currently fashionable categories which admit well-behaved smash
products; a coordinate-free approach suffices, as described in
\cite{el:ggs}, for example. The Pontryagin--Thom construction then
ensures that the double complex bordism and cobordism functors are
represented by the Thom spectrum $\MUMU$, which we label as $\DU$, and
the cobordism ring $\osdu$ is identified with the homotopy ring
$\pi_*(\DU)$. The transformation $\pi$ is induced by the product map on
$\MU$, whilst $\iota_\ell$ and $\iota_r$ are induced by the left and
right inclusion of $\MU$ in $\DU$ respectively, using the unit
$S^0\rightarrow\MU$ on the opposite factor.
 
We may also identify the homotopy ring of $\MUMU$ with the
$\osu$--algebra $\osu(\MU)$, adopting the convention of \cite{ad:shg}
(and most subsequent authors) in taking the argument as the second
factor. Of course, $\osu(\MU)$ is also the Hopf algebroid $A_*^U$ of
cooperations in complex bordism theory. The Thom isomorphism
$\osu(\MU)\cong\osu(\BU_+)$ provides a further description, whose ring
structure is induced by the Whitney sum map on the Grassmannian $\BU$;
it is commonly used to transfer the standard polynomial generator
$\beta_n$ in $\varOmega_{2n}^U(\BU)$ to the polynomial generator $b_n$
in $\varOmega_{2n}^U(\MU)$, for each $n\geq 0$. Monomials
$\beta^\omega$ in the $\beta_n$ are dual to the universal Chern
classes $c_\omega$ in $\ous(\BU)$, and monomials $b^\omega$ in the
$b_n$ are dual to the Landweber--Novikov operations $s_\omega$ in the
algebra of complex cobordism operations $A^*_U$. The {\it
Landweber--Novikov algebra} $S^*$ is the sub-Hopf algebra generated
additively by the $s_\omega$, with coproduct induced by the Cartan
formulae; its integral dual $S_*$ is the polynomial subalgebra
$\bZ[b_1,b_2,\dots]$ of $A_*^U$, with coproduct $\delta$ induced from
that of $A_*^U$ by restriction. We combine our isomorphisms as
\begin{equation}\label{canandh}
\osdu\cong\osu(\MU)\cong\osu(\BU_+), 
\end{equation} 
referring to the first as the {\it canonical isomorphism}, and to the
composition as $h$.  An analysis of the Pontryagin--Thom construction
confirms that $h$ maps the double cobordism class of any
$(M;\nu_\ell,\nu_r)$ to the cobordism class of the singular
$U$--manifold $\nu_r\colon M\rightarrow\BU$.

There are two complex orientation classes $x_\ell$ and $x_r$ in
$\varOmega^2_{DU}(\CPI)$, arising from the first Chern class $x$ in
$\varOmega^2_U(\CPI)$; indeed, $\DU$ is the universal example of a
{\it doubly complex oriented} spectrum. More generally, there are left
and right Chern classes $c_{\psi,\ell}$ and $c_{\omega,r}$ in
$\odus(\BU)$, dual to monomials $\beta^\psi_\ell$ and $\beta^\omega_r$
in $\osdu(\BU)$. We obtain mutually inverse expansions  
\begin{equation}\label{xrxs}
x_r=\sum_{n\geq 0}g_nx_\ell^{n+1}\spandsp x_\ell=\sum_{n\geq
0}\bar{g}_nx_r^{n+1} 
\end{equation}
in $\varOmega^2_{DU}(\CPI)$, where $g_n$ and $\bar{g}_n$ lie in
$\varOmega_{2n}^{DU}$ for all $n$ and are interchanged by the
involution $\chi$. For $n>0$ they are annihilated by the transformation
$\pi$, whilst $g_0=\bar{g}_0=1$. As documented in
\cite{bura:dcq}, the image of $g_n$ under the canonical isomorphism is
$b_n$, and the isomorphism $h$ of \eqref{canandh} therefore satisfies
$h(g_n)=\beta_n$ in $\varOmega^U_{2n}(\BU_+)$, for each $n\geq 0$.

These observations arise from minor manipulations with the
definitions, and suggest that we introduce the polynomial subalgebra
$G_*$ of $\osdu$, generated by the elements $g_n$ (or, equivalently,
by the elements $\bar{g}_n$) for $n\geq 0$. We may then
incorporate our previous remarks and formulate the geometric
viewpoint; we also appeal to \cite{ra:ocb}, recalling the construction
of singular manifolds $\beta\colon B^n\rightarrow\CPI$ to represent
$\beta_n$ in $\varOmega_{2n}^U(\CPI)$, where $B^n$ is an iterated
$2$--sphere bundle which admits a bounding $U$--structure for each
$n\geq 0$.  
\begin{prop}\label{giss} 
The canonical isomorphism identifies $G_*$ with the dual of the
Landweber--Novikov algebra $S_*$ in $A_*^U$; a representative for the
generator $g_n$ is given by $(B^n;\nu\oplus\beta^\perp,\beta)$, for
each $n\geq 0$. 
\end{prop}

We shall apply Proposition \ref{giss} to realize the coproduct and
antipode of $S_*$, given by  
\begin{equation}\label{lnstrucs}
\delta(b_n)=\sum_{k\geq 0}(b)_{n-k}^{n+1}\otimes b_k
\spandsp
\chi(b_n)=(b)^{-(n+1)}_n,
\end{equation}
and the left and right actions of $S^*$ on $S_*$, given
by 
\[
\langle y,s_\ell a\rangle=\langle\chi(s)y,a\rangle
\spandsp
\langle y,s_ra\rangle=\langle ys,a\rangle;
\]
here $s$ and $y$ lie in $S^*$, and the actions on $a$ in $S_*$
extend naturally to $A_*^U$. Equivalently, we may write
\begin{equation}\label{lractalt}
s_\ell a=\sum\langle \chi(s),a_1\rangle a_2
\spandsp
s_r a=\sum \langle s,a_2\rangle a_1
\end{equation}
where $\delta(a)=\sum a_1\otimes a_2$, confirming that either of the left
or right actions determines (and is determined by) the coproduct
$\delta$.  

We consider the algebra $A^*_{DU}\cong\odus(\DU)$ of operations in
double complex cobordism theory, whose continuous $\osdu$--dual is the
corresponding Hopf algebroid of cooperations
$A_*^{DU}\cong\osdu(\DU)$. An element $s$ of $S^*$ yields operations
$s_\ell\otimes 1$ and $1\otimes s_r$ by action on the first or second
factor $\MU$ of $\DU$, leading to the {\it left and right
Landweber--Novikov operations} $s_{\psi,\ell}\otimes 1$ and $1\otimes
s_{\omega,r}$, which commute in $A^*_{DU}$ by construction. It follows
that $A^*_{DU}$ contains the subalgebra $S^*_\ell\otimes S^*_r$, and
that $A_*^{DU}$ contains the subalgebra $S_{*,\ell}\otimes
S_{*,r}\cong\bZ[b_{j,\ell}\otimes 1,1\otimes b_{k,r}:j,k\geq 0]$;
these are integrally dual Hopf algebras. Of course $S^*_\ell\otimes
S^*_r$ acts on the coefficient ring $\osdu$, and we need only unravel
the definitions in order to express the result in terms of the
canonical isomorphism.
\begin{prop}\label{lislnrisr}
The canonical isomorphism identifies the actions of the algebras
$S^*_\ell\otimes 1$ and $1\otimes S^*_r$ on $\osdu$ with the left and
right actions of $S^*$ on $A_*^U$ respectively; in particular $G_*$ is
closed under the action of $S^*_\ell\otimes S^*_r$.
\end{prop}

Since $S^*$ is cocommutative, the image of the coproduct $\delta\colon
S^*\rightarrow S^*_\ell\otimes S^*_r$ defines a third subalgebra
$S^*_d$ of $A^*_{DU}$. The canonical isomorphism identifies the
resulting {\it diagonal action\/} of $S^*_d$ on $G_*$ with the {\it
adjoint action\/} of $S^*$ on $S_*$; this is fundamental to the
formation of the quantum double $\cD(S^*)$ \cite{ka:qg}, and underlies
the description of $\cD(S^*)$ as a subalgebra of $A^*_{DU}$
\cite{bu:smg}, \cite{bura:dcq}.

By analogy with standard cobordism theory the action of
$S^*_\ell\otimes S^*_r$ on $\osdu$ may be expressed in terms of
characteristic numbers, since the operation $s_{\psi,\ell}\otimes
s_{\omega,r}$ corresponds to the Chern class $c_{\psi,\ell}\otimes
c_{\omega,r}$ under the appropriate Thom isomorphism
$A^*_{DU}\cong\odus(\BU\times\BU_+)$. So the action of
$s_{\psi,\ell}\otimes s_{\omega,r}$ on the cobordism class of
$(M;\nu_\ell,\nu_r)$ is given by the Kronecker product
\begin{equation}\label{onptdoub}
\langle c_{\psi,\ell}(\nu_\ell)c_{\omega,r}(\nu_r),\sigma\rangle
\end{equation}
in $\osdu$, where $\sigma$ in $\osdu(M)$ is the canonical orientation
class represented by the identity map on $M$. The left and right
actions of $S^*$ are therefore given by restriction, yielding
$\langle c_{\psi,\ell}(\nu_\ell),\sigma\rangle$ and $\langle
c_{\omega,r}(\nu_r),\sigma\rangle$ respectively. Our procedure for
computing the actions of $S^*_\ell$ and $S^*_r$ on $G_*$ in Theorem
\ref{specss} is now revealed; we take the double $U$--cobordism class
of $(M;\nu_\ell,\nu_r)$, form the Poincar\'e duals of
$c_{\psi,\ell}(\nu_\ell)$ and $c_{\omega,r}(\nu_r)$ respectively, and
record the double $U$--cobordism classes of the resulting source
manifolds.  

\section{Bounded flag manifolds}\label{boflma}

In this section we introduce our family of bounded flag manifolds, and
discuss their topology in terms of a cellular calculus which is
intimately related to the Schubert calculus for classic flag manifolds.
Our description is couched in terms of nonsingular subvarieties,
anticipating applications to cobordism in the next section. We also
invest the bounded flag manifolds with certain canonical $U$-- and double
$U$--structures, and so relate them to our earlier constructions in
$\osdu$. Much of our notation differs considerably from that introduced
in \cite{ra:ocb}.

We shall follow combinatorial convention by writing $[n]$ for the set
of natural numbers $\{1,2,\dots,n\}$, equipped with the standard
linear ordering $<\,$. Every interval in the poset $[n]$ has the form
$[a,b]$ for some $1\leq a\leq b\leq n$, and consists of all $m$ satisfying
$a\leq m\leq b$; our convention therefore dictates that we abbreviate
$[1,b]$ to $[b]$. It is occasionally convenient to interpret $[0]$ as
the empty set, and $[\infty]$ as the natural numbers.  We work in the
context of the Boolean algebra $\mathcal{B}(n)$ of finite subsets of
$[n]$, ordered by inclusion. We decompose each such subset $Q\subseteq
[n]$ into maximal subintervals $I(1)\cup\dots\cup I(s)$, where
$I(j)=[a(j),b(j)]$ for $1\leq j\leq s$, and assign to $Q$ the monomial
$b^\omega$, where $\omega_i$ records the number of intervals $I(j)$ of
cardinality $i$ for each $1\leq i\leq n$; we refer to $\omega$ as
the {\it type} of $Q$, noting that it is independent of the choice of
$n$. We display the elements of $Q$ in increasing order as $\{q_i:1\leq
i\leq d\}$, and abbreviate the complement $[n]\setminus Q$ to $Q'$. We
also write $I(j)^+$ for the subinterval $[a(j),b(j)+1]$ of $[n+1]$, and
$Q^\wedge$ for $Q\cup\{n+1\}$. It is occasionally convenient to set $b(0)$
to $0$ and $a(s+1)$ to $n+1$. 

We begin by recalling standard constructions of complex flag
manifolds and some of their simple properties, for which a helpful
reference is \cite{hi:gcg}. We work in an ambient complex inner
product space $Z_{n+1}$, which we assume to be invested with a 
preferred orthonormal basis $z_1$,\dots, $z_{n+1}$, and
we write $Z_E$ for the subspace spanned by the vectors $\{z_e:e\in
E\}$, where $E\subseteq [n+1]$. We abbreviate $Z_{[a,b]}$ to
$Z_{a,b}$ (and $Z_{[b]}$ to $Z_b$) for each $1\leq a<b\leq n+1$, and   
write $\CP(Z_E)$ for the projective space of lines in $Z_E$.
We let $V-U$ denote the orthogonal complement of $U$ in $V$ for any 
subspaces $U<V$ of $Z_{n+1}$, and we regularly abuse notation by 
writing 0 for the subspace which consists only of the zero vector. A
complete flag $V$ in $Z_{n+1}$ is a sequence of proper subspaces   
\[
0=V_0<V_1<\dots<V_i<\dots<V_{n}<V_{n+1}=Z_{n+1},
\]
of which the {\it standard\/} flag 
$Z_0<\dots<Z_i<\dots<Z_{n+1}$ is a specific example. The flag
manifold $F(Z_{n+1})$ is the set of all flags in $Z_{n+1}$, 
topologized as the quotient $U(n+1)/T$ of the unitary group $U(n+1)$
by its maximal torus.

The flag manifold is a nonsingular complex projective algebraic
variety of dimension $\tbinom{n+1}{2}$, whose cells $e_\alpha$ are
even dimensional, indexed by elements of the symmetric group
$\mathfrak{S}_{n+1}$, and partially ordered by the decomposition of
$\alpha$ into a product of transpositions. The closure of every
$e_\alpha$ is an algebraic subvariety, generally singular, known as
the {\it Schubert variety} $X_\alpha$. Whether considered as cells or
subvarieties, the $e_\alpha$ define a basis for the integral homology
and cohomology groups $H_*(F(Z_{n+1}))$ and $H^*(F(Z_{n+1}))$, which
are integrally dual. The manipulation of cup and cap products and
Poincar\'e duality in these terms is known as the {\it Schubert
calculus} for $F(Z_{n+1})$.

An alternative description of $H^*(F(Z_{n+1}))$ is provided by Borel's
computations with the {\it characteristic homomorphism}
$H^*(BT)\rightarrow H^*(U(n+1)/T)$, induced by the canonical torus
bundle $U(n+1)/T\rightarrow BT$. Noting that $H^*(BT)$ is a polynomial
algebra on two dimensional generators $x_i$ for $1\leq i\leq n$, Borel
identifies $H^*(F(Z_{n+1}))$ with the ring of {\it coinvariants} under
the action of $\mathfrak{S}_{n+1}$. In this context, $x_i$ is the
first Chern class of the line bundle over $F(Z_{n+1})$ obtained by
associating $V_i-V_{i-1}$ to each flag $V$.

The interaction between the Schubert and Borel descriptions of the
cohomology of $F(Z_{n+1})$ is a fascinating area of combinatorial algebra
and has led to a burgeoning literature on the subject of {\it Schubert
polynomials}, beautifully introduced in MacDonald's book \cite{ma:nsp}. 

We call a flag $U$ in $Z_{n+1}$ {\it bounded\/} if each
$i$--dimensional component $U_i$ contains the first $i-1$ basis vectors
$z_1$,\dots, $z_{i-1}$, or equivalently, if $Z_{i-1}<U_i$ for every
$1\leq i\leq n+1$. We define the {\it bounded flag manifold}
$B(Z_{n+1})$ to be the set of all bounded flags in $Z_{n+1}$,
topologized as a subvariety of $F(Z_{n+1})$; it is straightforward to
check that $B(Z_{n+1})$ is nonsingular, and has dimension $n$. Clearly
$B(Z_2)$ is isomorphic to the projective line $\CP(Z_2)$ with the
standard complex structure, whilst $B(Z_1)$ consists solely of the
trivial flag. We occasionally abbreviate $B(Z_{n+1})$ to $B_n$, in
recognition of its dimension. 

The algebraic torus $(C^*)^n$ is contained in $Z_n$, and each of
its points $t$ determines a line $L_t<Z_{n+1}$ with basis vector
$t+z_{n+1}$. We may therefore embed $(C^*)^n$ in $B(Z_{n+1})$ as an
open dense subset, by assigning the bounded flag 
\[
0<L_t<L_t\oplus Z_1<\dots<L_t\oplus Z_i<\dots<L_t\oplus
Z_{n-1}<Z_{n+1}
\]
to each $t$. The standard action of $(C^*)^n$ on this torus extends to
the whole of $B(Z_{n+1})$ by coordinatewise multiplication on $Z_n$
(fixing $z_{n+1}$), and therefore imposes a canonical toric variety
structure \cite{fu:itv}. 

There is a map $p_h\colon B(Z_{n+1})\rightarrow B(Z_{h+1,n+1})$
for each $1\leq h\leq n$, defined by factoring out $Z_h$. Thus 
$p_h(U)$ is given by
\[
0<U_{h+1}-Z_h<\dots<U_i-Z_h<\dots<U_n-Z_h<Z_{h+1,n+1}
\]
for each bounded flag $U$. Since $Z_{i-1}<U_i$ for all 
$1\leq i\leq n+1$, we deduce that
$Z_{h+1,i-1}<U_i-Z_h$ for all $i>h+1$, ensuring that
$p_h(U)$ is indeed bounded. We may readily check that $p_h$ is
the projection of a fiber bundle, with fiber $B(Z_{h+1})$.
In particular, $p_1$ has fiber $\CP(Z_2)$, and after $n-1$
applications we may exhibit $B(Z_{n+1})$ as an iterated bundle 
\begin{equation}\label{iterbdle}
B(Z_{n+1})\rightarrow\cdots\rightarrow B(Z_{h,n+1})
\rightarrow\cdots\rightarrow B(Z_{n,n+1}),
\end{equation}
where the fiber of each map is isomorphic to $\CP^1$. This
construction was introduced in \cite{ra:ocb}.

We define maps $q_h$ and $r_h\colon B(Z_{n+1})\rightarrow
\CP(Z_{h,n+1})$ by letting $q_h(U)$ and $r_h(U)$ be the respective
lines $U_h-Z_{h-1}$ and $U_{h+1}-U_h$, for each $1\leq h\leq n$. We
remark that $q_h=q_1\cdot p_{h-1}$ and $r_h=r_1\cdot p_{h-1}$ for all
$h$, and that the appropriate $q_h$ and $r_h$ may be assembled into
maps $q_Q$ and $r_Q\colon
B(Z_{n+1})\rightarrow\times_Q\CP(Z_{h,n+1})$, where $h$ varies over an
arbitrary subset $Q$ of $[n]$. In particular, $q_{[n]}$ is an
embedding which associates to each flag $U$ the $n$--tuple
$(U_1,\dots,U_h-Z_{h-1},\dots,U_n-Z_{n-1})$, and describes
$B(Z_{n+1})$ as a projective algebraic variety.

We proceed by analogy with the Schubert calculus for $F(Z_{n+1})$. To
every flag $U$ in $B(Z_{n+1})$ we assign the {\it support} $S(U)$, given by 
$\{j\in [n]:U_j\neq Z_j\}$, and consider the subspace  
\[
e_Q=\{U\in B(Z_{n+1}):S(U)=Q\}
\]
for each $Q$ in the Boolean algebra $\mathcal{B}(n)$. For example,
$e_\emptyset$ is the singleton consisting of the standard flag.
\begin{lem}\label{opcell} 
The subspace $e_Q\subset B(Z_{n+1})$ is an open cell of dimension
$2|Q|$, whose closure $X_Q$ is the union of all $e_R$ for which
$R\subseteq Q$ in $\mathcal{B}(n)$.     
\end{lem}
\begin{proof}
If $Q=\cup_jI(j)$, then $e_Q$ is
homeomorphic to the cartesian product $\times_je_{I(j)}$, so it
suffices to assume that $Q$ is an interval $[a,b]$. If $U$ lies in
$e_{[a,b]}$ then $U_{a-1}=Z_{a-1}$ and $U_{b+1}=Z_{b+1}$ certainly
both hold;
thus $e_{[a,b]}$ consists of those flags $U$ for which $q_j(U)$ is a
fixed line $L$ in $\CP(Z_{a,b+1})\setminus\CP(Z_{a,b})$ for all $a\leq 
j\leq b$. Therefore $e_{[a,b]}$ is a $2(b-a+1)$--cell, as sought.
Obviously $e_R\subset X_{[a,b]}$ for each $R\subseteq Q$, so it
remains only to observe that the limit of a sequence of flags in
$e_{[a,b]}$ cannot have fewer components satisfying $U_j=Z_j$, and
must therefore lie in $e_R$ for some $R\subseteq [a,b]$. 
\end{proof}
Clearly $X_{[n]}$ is $B(Z_{n+1})$, so that Lemma \ref{opcell} provides a 
CW decomposition for $B_n$ with $2^n$ cells. 

We now prove that all the subvarieties $X_Q$ are nonsingular, in
contrast to the situation for $F(Z_{n+1})$. 

\begin{prop}\label{xqsnonsin}
For any $Q\subseteq [n]$, the subvariety $X_Q$ is diffeomorphic to the
cartesian product $\times_jB(Z_{I(j)^+})$.
\end{prop}
\begin{proof}
We may define a smooth embedding
$i_Q\colon\times_jB(Z_{I(j)^+})\rightarrow B(Z_{n+1})$ by choosing the
components of $i_Q(U(1),\dots,U(s))$ to be
\begin{equation}\label{xqhelp}
T_k=
\begin{cases}
Z_{a(j)-1}\oplus U(j)_i&\text{if $k=a(j)+i-1$ in $I(j)$}\\
Z_k&\text{if $k\in [n+1]\setminus Q$},                           
\end{cases}
\end{equation}
where $U(j)_i<Z_{I(j)^+}$ for each $1\leq i\leq b(j)-a(j)+1$; 
the resulting flag is indeed bounded, since 
$Z_{a(j),a(j)+i-1}<U(j)_i$ holds for all such $i$ and $1\leq j\leq s$.
Any flag $T$ in $B(Z_{n+1})$ for which $S(T)\subseteq Q$ must be of
the form \eqref{xqhelp}, so that $i_Q$ has image $X_Q$, as required. 
\end{proof}

We may therefore interpret the set 
\[
\mathcal{X}(n)=\{X_Q:Q\in\mathcal{B}(n)\}
\]
as a Boolean algebra of nonsingular subvarieties of $B(Z_{n+1})$,
ordered by inclusion, on which the support function
$S\colon\mathcal{X}(n)\rightarrow\mathcal{B}(n)$ induces an isomorphism
of Boolean algebras. Moreover, whenever $Q$ has type $\omega$ then $X_Q$
is isomorphic to the cartesian product   
$B_1^{\omega_1}B_2^{\omega_2}\dots B_n^{\omega_n}$, and so may be
abbreviated to $B^\omega$. In this important sense, $S$ preserves types.
We note that the complex dimension $|Q|$ of $X_Q$ may be written as
$|\omega|$.    

The following quartet of lemmas is central to our computations in
\S\ref{nostdu}.  
\begin{lem}\label{xtrans}
The map $r_{Q'}\colon
B(Z_{n+1})\rightarrow\times_{Q'}\CP(Z_{h,n+1})$ is
transverse to the subvariety $\times_{Q'}\CP(Z_{h+1,n+1})$,
whose inverse image is $X_Q$.
\end{lem}
\begin{proof}
Let $T$ be a flag in $B(Z_{n+1})$. Then $r_h(T)$ lies $\CP(Z_{h+1,n+1})$
if and only if $T_{h+1}=T_h\oplus L_h$ for some line $L_h$ in
$Z_{h+1,n+1}$. Since $Z_h<T_{h+1}$, this condition is equivalent to
requiring that $T_h=Z_h$, and the proof is completed by allowing $h$ to
range over $Q'$. 
\end{proof}
\begin{lem}\label{ytrans}
The map $q_{Q'}\colon
B(Z_{n+1})\rightarrow\times_{Q'}\CP(Z_{h,n+1})$ is
transverse to the subvariety $\times_{Q'}\CP(Z_{h+1,n+1})$,
whose inverse image is diffeomorphic to $B(Z_{Q^\wedge})$.
\end{lem}
\begin{proof}
Let $T$ be a flag in $B(Z_{n+1})$ such that $q_h(T)$ lies
$\CP(Z_{h+1,n+1})$, which occurs if and only if $T_h=Z_{h-1}\oplus L_h$ 
for some line $L_h$ in $Z_{h+1,n+1}$. Whenever this equation holds for
all $h$ in some interval $[a,b]$, we deduce that $L_h$ actually lies
in $Z_{b+1,n+1}$. Thus we may describe $T$ globally by
\[
T_k=Z_{[k-1]\setminus Q}\oplus U_i,
\]
where $U_i$ lies in $Z_{Q^\wedge}$, and $i$ is $k-|[k-1]\setminus Q|$.
Clearly $U_{i-1}<U_i$ and $Z_{\{q_1,\dots,q_{i-1}\}}<U_i$ for all
appropriate $i$, so that $U$ lies in $B(Z_{Q^\wedge})$. We may now identify 
the required inverse image with the image of the natural smooth
embedding $j_Q\colon B(Z_{Q^\wedge})\rightarrow B(Z_{n+1})$, as sought.
\end{proof}

We therefore define $Y_Q$ to consist of all flags $T$ for which the
line $T_h-Z_{h-1}$ lies in $Z_{Q^\wedge}$ for every $h$ in $Q'$. It
follows that $Y_Q$ is isomorphic to $B_k$ whenever $Q$ has cardinality
$k$; for example, $Y_{[n]}$ is $B(Z_{n+1})$ itself and $Y_\emptyset$
consists of the single flag determined by $T_1=Z_{n+1}$. The set
\[
\mathcal{Y}(n)=\{Y_Q:Q\in\mathcal{B}(n)\}
\]
is also a Boolean algebra of nonsingular subvarieties. 

\begin{lem}\label{powertr}
For any $1\leq m\leq n-h$, the map 
$q_h\colon B(Z_{n+1})\rightarrow\CP(Z_{h,n+1})$ is
transverse to the subvariety $\CP(Z_{h+m,n+1})$,
whose inverse image is diffeomorphic to $Y_{[h,h+m-1]'}$.
\end{lem}
\begin{proof}
Let $T$ be a flag in $B(Z_{n+1})$ such that $q_h(T)$ lies
$\CP(Z_{h+m,n+1})$, which occurs if and only if $T_h=Z_{h-1}\oplus L_h$ 
for some line $L_h$ in $Z_{h+m,n+1}$. Following the proof of Lemma
\ref{ytrans} we immediately identify the required inverse image
with $Y_{[h-1]\cup [h+m,n]}$, as sought.
\end{proof}

\begin{lem}\label{intersecs}
The following intersections in $B(Z_{n+1})$ are transverse:
\begin{gather*}
X_Q\cap X_R=X_{Q\cap R}\sands Y_Q\cap Y_R=Y_{Q\cap R}
\sts{whenever}Q\cup R=[n],\\
\quad\text{and}\quad X_Q\cap Y_{R}=
\begin{cases}
X_{Q,R}&\text{if $Q\cup R=[n]$}\\
\;\emptyset&\text{otherwise},
\end{cases}
\end{gather*}
where $X_{Q,R}$ denotes the submanifold $X_{Q\cap R}\subseteq
B(Z_{R^\wedge})$. 
Moreover, $m$ copies of $Y_{\{h\}'}$ may be made self-transverse so that 
\[
Y_{\{h\}'}\cap\dots\cap Y_{\{h\}'}=Y_{[h,h+m-1]'}
\]
for each $1\leq h\leq n$ and $1\leq m\leq n-h$.
\end{lem}
\begin{proof}
The first three formulae follow directly from the definitions, and
dimensional considerations ensure that the intersections are transverse.
The manifold $X_{Q,R}$ is diffeomorphic to
$\dtimes_jY_{R(j)}$ as a submanifold of $B(Z_{n+1})$, where
$Q=\cup_jI(j)$ and $R(j)=I(j)\cap R$ for each $1\leq j\leq s$.

Since $Y_{\{h\}'}$ is defined by the single constraint
$U_h=Z_{h-1}\oplus L_h$, where $L_h$ is a line in $Z_{h+1,n+1}$, 
we may deform the embedding $j_{\{h\}'}$ (through smooth embeddings, in 
fact) to $m-1$ further embeddings in which the $L_h$ is constrained to
lie in $Z_{[h,n+1]\setminus\{h+i-1\}}$, for each $2\leq i\leq m$. The
intersection of the $m$ resulting images is determined by the single
constraint $L_h<Z_{h+m,n+1}$, and the result follows by applying Lemma
\ref{powertr}.
\end{proof}

It is illuminating to consider the toric structure of $B(Z_{n+1})$ in
these terms.
\begin{prop}\label{toptor}\hfill
\begin{enumerate}
\item[\rm (1)]
For each $1\leq h\leq n$, the projection $q_h\colon
B(Z_{n+1})\rightarrow\CP(Z_{h,n+1})$ is equivariant with respect to an
action of the torus $(C^*)^{n-h+1}$, and the equivariant filtration 
\[
\CP(Z_h)\subset\dots\subset\CP(Z_{h,i})\subset\dots\subset\CP(Z_{h,n+1})
\]
lifts to an equivariant filtration of the irregular values of $q_h$.
\item[\rm (2)]
The quotient of $B(Z_{n+1})$ by the action of the compact torus $T^n$
is homeomorphic to the $n$--cube $I^n$.
\end{enumerate}
\end{prop}
\begin{proof}
For (1), we choose the subtorus of $(C^*)^n$ in which the first $h-1$
coordinates are $1$; in particular, when $h=1$ the result refers to
toric structures on $B(Z_{n+1})$ and $\CP^n$. 
For (2), we proceed inductively from the observation that the
invariant submanifolds of the action of $T^n$ are the subvarieties
$X_{Q\setminus R,Q}$ for all pairs $R\subseteq Q\subseteq [n]$; in
particular, the fixed points are standard flags in the subvarieties 
$B(Z_{Q^\wedge})$, and so display the vertices of the quotient in
bijective correspondence with the subsets $Q$.
\end{proof}

The second part of Proposition \ref{toptor} refers to the structure of
$B(Z_{n+1})$ as a toric {\it manifold\/} \cite{daja:cpc}, and may be
extended by algebraic geometers to a more detailed description of the
associated fan \cite{fu:itv}. 

\section{Normal structures and duality}\label{nostdu}

In this section we describe the basic $U$-- and double $U$--structures on
the varieties $X_Q$, and compute their cobordism rings. We pay special
attention to Poincar\'e duality, which makes delicate use of the normal
structures and is of central importance to our subsequent applications. 

We consider complex line bundles $\gamma_i$ and $\rho_i$
over $B(Z_{n+1})$, classified respectively by the maps $q_i$
and $r_i$ for each $1\leq i\leq n$. We set $\gamma_0$ to $0$ and
$\rho_0$ to $\gamma_1$, which are compatible with the choices above and
enable us to write 
\begin{equation}\label{bunprops}
\gamma_i\oplus\rho_i\oplus\rho_{i+1}\oplus\dots\oplus\rho_n
\cong\bC^{\,n-i+2}
\end{equation} 
for every $0\leq i\leq n$. We may follow \cite{ra:ocb} in using
\eqref{iterbdle} to obtain an expression of the form
$\tau\oplus\bR\cong(\oplus_{i=2}^{n+1}\gamma_i)\oplus\bR$ for the
tangent bundle of $B(Z_{n+1})$, as prophesied by the toric structure;
so \eqref{bunprops} leads to an isomorphism
$\nu\cong\oplus_{i=2}^n(i-1)\rho_i$. We refer to the resulting
$U$--structure as the {\it basic\/} $U$--structure on $B(Z_{n+1})$. We
emphasize that these isomorphisms are of real bundles only, and
that the basic $U$--structure is not compatible with any complex 
structure on the underlying variety. On $B(Z_2)$, for example, the
basic $U$--structure is that of a 2--sphere $S^2$, rather than
$\CP^1$. Indeed, the basic $U$--structure on $B(Z_{n+1})$ extends over
the 3--disc bundle associated to $\gamma_1\oplus\bR$ for all values of
$n$, so that $B(Z_{n+1})$ represents zero in $\varOmega^U_{2n}$.   

By virtue of \eqref{bunprops} we may introduce the double $U$--structure 
$(\bigoplus_{i=1}^ni\rho_i,\gamma_1)$, which we again label {\it
basic}; equivalently, we rewrite $\nu_\ell$ as 
$\gamma=-(\gamma_1\oplus\dots\oplus\gamma_n)$. The basic double 
$U$--structure does not bound, however, as we shall see in Proposition
\ref{bnisbn}. Given any cartesian product of manifolds $B(Z_{n+1})$,
we also refer to the product of basic structures as basic.  
\begin{prop}\label{bnisbn}
With the basic double $U$--structure, $B(Z_{n+1})$ represents $g_n$ in
$\varOmega^{DU}_*$; if $\nu_\ell$ and $\nu_r$ are interchanged, it
represents $\bar{g}_n$.   
\end{prop}
\begin{proof}
It suffices to apply Proposition \ref{giss} for $g_n$, because the
bundle $\gamma_1$ over $B(Z_{n+1})$ coincides with the bundle $\beta$
of \cite{ra:ocb} over $B^n$. The result for $\bar{g}_n$ follows
by applying the involution $\chi$. 
\end{proof}
\begin{cor}\label{xqgenln}
The cobordism classes of the basic double $U$--manifolds $X_Q$ give an
additive basis for $G_*$ as $Q$ ranges over finite subsets of $[\infty]$.  
\end{cor}
\begin{proof}
It suffice to combine Propositions \ref{xqsnonsin} and \ref{bnisbn},
remarking that $X_Q$ represents $g^\omega$ whenever $Q$ has type
$\omega$.  
\end{proof}
Henceforth we shall insist that $B_n$ denotes $B(Z_{n+1})$ (or any
isomorph) equipped exclusively with the basic double $U$--structure.

\begin{prop}\label{xnystrucs}
Both $\mathcal{X}(n)$ and $\mathcal{Y}(n)$ are Boolean algebras of
basic $U$--submanifolds, in which the intersection formulae of Lemma
{\rm \ref{intersecs}} respect the basic $U$--structures. 
\end{prop}
\begin{proof}
It suffices to prove that the pullbacks in Lemmas \ref{xtrans},
\ref{ytrans} and \ref{powertr} are compatible with the basic
$U$--structures. Beginning with Lemma \ref{xtrans}, we note that
whenever $\rho_h$ over $B(Z_{n+1})$ is restricted by $i_Q$ 
to a factor $B(Z_{I(j)^+})$, we
obtain $\rho_{k+1}$ if $h=a(j)+k$ lies in $I(j)$ and $\gamma_1$ if
$h=a(j)-1$; for all other values of $h$, the
restriction is trivial. Since the construction of Lemma \ref{xtrans}
identifies $\nu(i_Q)$ with the restriction of $\oplus_h\rho_h$ as $h$
ranges over $Q'$, we infer an isomorphism
$\nu(i_Q)\cong(\times_j\gamma_1)\oplus\bC^{\,n-j-|Q|}$ over $X_Q$
(unless $1\in Q$, in which case the first $\gamma_1$ is trivial).
Appealing to \eqref{bunprops}, we then verify that this is compatible
with the basic structures in the isomorphism 
$\nu^{X_Q}\cong i_Q^*(\nu^{B(Z_{n+1})})\oplus\nu(i_Q)$, 
as claimed. The proofs for Lemmas \ref{ytrans} and \ref{powertr} are
similar, noting that the 
restriction of $\rho_h$ to $Y_Q$ is $\rho_k$ if $h=q_k$ lies in
$Q$, and is trivial otherwise, and that the restriction of $\gamma_h$ is
$\gamma_k$ if $h=q_k$ lies in $Q$, and is $\gamma_{k+1}$ if $q_k$ is the
greatest element of $Q$ for which $h>q_k$ (meaning $\gamma_1$ if
$h<q_1$, and the trivial bundle if $h>q_k$ for all $k$). Since the
construction of Lemma \ref{ytrans} identifies  
$\nu(j_Q)$ with the restriction of $\oplus_h\gamma_h$ as $h$ ranges over
$Q'$, we infer an isomorphism 
\begin{equation}\label{nujq}
\nu(j_Q)\cong\bigoplus_{j=1}^{s+1}\Big(a(j)-b(j-1)-1\Big)\gamma_{c(j)}
\end{equation}
over $Y_Q$, where $c(j)=j+\sum_{i=0}^{j-1}(b(i)-a(i))$. This isomorphism
is also compatible with the basic structures in $\nu^{Y_Q}\cong
j_Q^*(\nu^{B(Z_{n+1})})\oplus\nu(j_Q)$, once more by appeal to
\eqref{bunprops}. 
\end{proof}
The corresponding results for double $U$--structures are more subtle,
since we are free to choose our splitting of $\nu(i_Q)$ and $\nu(j_Q)$
into left and right components. 
\begin{cor}\label{xnydub}
The same results hold for double $U$--structures with respect to the
splittings  
$\nu(i_Q)_\ell=0$ and $\nu(i_Q)_r=\nu(i_Q)$, and
$\nu(j_Q)_\ell=\nu(j_Q)$ and $\nu(j_Q)_r=0$.  
\end{cor}
\begin{proof}
One extra fact is required in the calculation for $i_Q$, namely that 
$\gamma_1$ on $B(Z_{n+1})$ restricts trivially to $X_Q$ (or to
$\gamma_1$ if $1\in Q$).  
\end{proof}

At this juncture we may identify the inclusions of $X_Q$ in
$F(Z_{n+1})$ with certain of the desingularizations introduced by Bott
and Samelson \cite{bosa:atm}. For example, $X_{[n]}$ is the 
desingularization of the Schubert variety
$X_{(n+1,1,2,\dots,n)}$, and the resolution map is actually an
isomorphism in this case. Moreover, the corresponding $U$--cobordism
classes form the cornerstone of Bressler and Evens's calculus for
$\ous(F(Z_{n+1}))$. In both of these applications, however, the
underlying complex manifold structures suffice. The basic
$U$--structures become vital when investigating the Landweber--Novikov
algebra (and could also have been used in \cite{brev:scc}, although an
alternative calculus would result). We leave the details to interested
readers.

We now use the basic structures on $X_Q$ to investigate Poincar\'e
duality in bordism and cobordism, beginning with the CW decomposition
for $B(Z_{n+1})$ which stems from Lemma \ref{opcell}. Since the cells
$e_Q$ occur only in even dimensions, the corresponding homology
classes $x_Q^H$ form a basis for the integral homology groups
$H_*(B(Z_{n+1}))$ as $Q$ ranges over $\mathcal{B}(n)$. Applying
$\Hom_{\bZ}$ determines a dual basis $\hd(x_Q^H)$ for the cohomology
$H^*(B(Z_{n+1}))$; we delay clarifying the cup product structure until
after Theorem \ref{cobordz} below, although it may also be deduced
directly from the toric properties of $B(Z_{n+1})$.

We introduce the complex bordism classes $x_Q$ and $y_Q$ in
$\varOmega^U_{2|Q|}(B(Z_{n+1}))$, represented respectively by the
inclusions $i_Q$ and $j_Q$ of the subvarieties $X_Q$ and $Y_Q$ with
their basic $U$--structures. By construction, the fundamental class in    
$H_{2|Q|}(X_Q)$ maps to $x_Q^H$ in $H_{2|Q|}(B(Z_{n+1}))$ under
$i_Q$; thus $x_Q$ maps to $x_Q^H$ under the Thom homomorphism 
$\osu(B(Z_{n+1}))\rightarrow H_*(B(Z_{n+1}))$. The \ahss\ for
$\osu(B(Z_{n+1}))$ therefore collapses, and the classes $x_Q$ form an
$\osu$--basis as $Q$ ranges over $\mathcal{B}(n)$. The classes
$x_{[n]}$ and $y_{[n]}$ coincide, since they are both represented by
the identity map. They constitute the {\it basic fundamental class} in 
$\varOmega^U_{2n}(B(Z_{n+1}))$, with respect to which the Poincar\'e 
duality isomorphism is given by 
\[
\pd(w)=w\cap x_{[n]}
\]
in $\varOmega^U_{2(n-d)}(B(Z_{n+1}))$, for any $w$ in
$\varOmega_U^{2d}(B(Z_{n+1}))$. 

An alternative source of elements in $\varOmega_U^2(B(Z_{n+1}))$ is
provided by the Chern classes 
\[
x_i=c_1(\gamma_i)\spandsp y_i=c_1(\rho_i)
\]
for each $1\leq i\leq n$. It follows from \eqref{bunprops} that
\begin{equation}\label{xitoy}  
x_i=-y_i-y_{i+1}-\cdots-y_n
\end{equation}
for every $i$. Given $Q\subseteq [n]$, we write $\prod_Qx_h$ as $x^Q$
and $\prod_Qy_h$ as $y^Q$ in $\varOmega_U^{2|Q|}(B(Z_{n+1}))$, where
$h$ ranges over $Q$ in both products. 

We may now discuss the implications of our intersection results of
Lemma \ref{intersecs} for the structure of $\varOmega_U^*(B(Z_{n+1}))$.
It is convenient (but by no means necessary) to use Quillen's
geometrical interpretation of cobordism classes, which provides a
particularly succinct description of cup and cap products and
Poincar\'e duality, and is conveniently summarized in \cite{brev:scc}. 
\begin{thm}\label{cobordz}
The complex bordism and cobordism of $B(Z_{n+1})$ satisfy
\begin{enumerate}
\item[(1)]
$\pd(x^{Q'})=y_Q$ and\/ $\pd(y^{Q'})=x_Q$;
\item[(2)]
the elements $\{y_Q:Q\subseteq [n]\}$ form an $\osu$--basis for
$\osu(B(Z_{n+1}))$;
\item[(3)]
$\hd(x_Q)=x^Q$ and\/ $\hd(y_Q)=y^Q$;
\item[(4)]
there is an isomorphism of rings
\[
\varOmega_U^*(B(Z_{n+1}))\cong\osu[x_1,\dots,x_n]/(x_i^2=x_ix_{i+1}),
\]
where $i$ ranges over $[n]$ and $x_{n+1}$ is interpreted as $0$.
\end{enumerate}
\end{thm}    
\begin{proof}
For (1), we apply Lemma \ref{ytrans} and Proposition \ref{xnystrucs} 
to deduce that $x^{Q'}$ in $\varOmega_U^{2|Q'|}(B(Z_{n+1}))$ is the
pullback of the Thom class under the collapse map onto $M(\nu(j_Q))$. 
Hence $x^{Q'}$ is represented geometrically by the inclusion
$j_Q\colon Y_Q\rightarrow B(Z_{n+1})$, and therefore $\pd(x^{Q'})$ is
represented by the same singular $U$--manifold in
$\varOmega_{2|Q|}(B(Z_{n+1}))$. Thus $\pd(x^{Q'})=y_Q$. An identical
method works for $\pd(y^{Q'})$, by applying Lemma \ref{xtrans}.
For (2), we have already shown that the $x_Q$ form an $\osu$--basis for
$\osu(B(Z_{n+1}))$. Thus by (1) the $y^Q$ form a basis for
$\ous(B(Z_{n+1}))$, and therefore so do the $x^Q$ by \eqref{xitoy}; the
proof is concluded by appealing to (1) once more.
To establish (3), we remark that the cap product $x^Q\cap x_R$ is
represented geometrically by the fiber product of $j_{Q'}$ and $i_R$,
and is therefore computed by the intersection theory of Lemma
\ref{intersecs}. Bearing in mind the crucial fact that each basic
$U$--structure bounds (except in dimension zero!), we obtain
\begin{equation}\label{unotdu}
\langle x^Q,x_R\rangle =\delta_{Q,R}
\end{equation}
and therefore that $\hd(x_Q)=x^Q$, as sought. The result for
$\hd(y_Q)$ follows similarly.
To prove (4) we note that it suffices to obtain the product formula
$x_i^2=x_ix_{i+1}$, since we have already demonstrated that the
monomials $x^Q$ form a basis in (2). Now $x_i$ and $x_{i+1}$ are
represented geometrically by $Y_{\{i\}'}$ and $Y_{\{i+1\}'}$
respectively, and products are represented by intersections; according
to Lemma \ref{intersecs} (with $m=2$), both $x_i^2$ and $x_ix_{i+1}$ are
therefore represented by the same subvariety $Y_{\{i,i+1\}'}$, so long
as $1\leq i<n$. When $i=n$ we note that $x_n$ pulls back from $\CP^1$,
so that $x_n^2=0$, as required.
\end{proof}

For any $Q\subseteq [n]$, we obtain the corresponding
structures for the complex bordism and cobordism of $X_Q$ by applying
the K\"unneth formula to Theorem \ref{cobordz}. Using the same
notation as in $B(Z_{n+1})$ for any cohomology class which restricts
along (or homology class which factors through) the inclusion $i_Q$,
we deduce, for example, a ring isomorphism     
\begin{equation}\label{cobordxq}
\varOmega_U^*(X_Q)\cong\osu[x_i:i\in Q]/(x_i^2=x_ix_{i+1}),
\end{equation}
where $x_i$ is interpreted as $0$ for all $i\notin Q$.

The relationship between the classes $x_i$ and $y_i$ in
$\varOmega_U^*(B(Z_{n+1}))$ is described by \eqref{xitoy}, but may be
established directly by appeal to the third formula of Lemma
\ref{intersecs}, as in the proof of Theorem \ref{cobordz}; for
example, we deduce immediately that $x_iy_i=0$ for all $1\leq i\leq
n$. When applied with arbitrary $m$, the fourth formula of Lemma
\ref{intersecs} simply iterates the quadratic relations, and produces
nothing new. 

The results of Theorem \ref{cobordz} extend to any complex oriented
cohomology theory as usual; in particular, we may substitute double
complex cobordism, so long as we choose left or right Chern classes
consistently throughout. To understand duality, however, we must also
attend to the choice of splittings provided by Corollary \ref{xnydub},
and the failure of formulae such as \eqref{unotdu} because the
manifolds $B_n$ are no longer double $U$--boundaries. Since, by
\eqref{onptdoub}, duality lies at the heart of our applications to the
Landweber--Novikov algebra, we treat these issues with care below.

We are particularly interested in the left and right Chern classes 
$x^Q_\ell$, $y^Q_\ell$, $x^Q_r$ and $y^Q_r$ in
$\varOmega^{2|Q|}_{DU}(B_n)$, and we seek economical geometric
descriptions of their Poincar\'e duals. We continue to write
$x_R$ and $y_R$ in $\varOmega_{2|R|}^{DU}(B_n)$ for the homology
classes represented by the respective inclusions of $X_R$ and $Y_R$ with 
their basic double $U$--structures.  
\begin{prop}\label{pdinbns}
In $\varOmega_{2(n-|Q|)}^{DU}(B_n)$, we have that
\[
\pd(x^{Q'}_\ell)=y_Q\spandsp \pd(y^{Q'}_r)=x_Q,
\]
whilst $\pd(x^{Q'}_r)$ and $\pd(y^{Q'}_\ell)$ are represented by the
inclusion of $Y_Q$ and $X_Q$ with the respective double $U$--structures  
\[
\big(\nu^{Y_Q}-(\nu(j_Q)\oplus\gamma_1),\;
\nu(j_Q)\oplus\gamma_1\big)
\sands
(\nu^{X_Q}-i^*_Q\gamma_1,\; i^*_Q\gamma_1),
\]
for all $n\geq 0$.
\end{prop} 
\begin{proof}
The first two formulae follow at once from Corollary \ref{xnydub},
by analogy with (1) of Theorem \ref{cobordz}. The second
require the interchange of the left and right components of
the normal bundles of $j_Q$ and $i_Q$ respectively, plus the
observation that $j_Q^*(\gamma_1)$ is always $\gamma_1$, whatever $Q$.
\end{proof}

Proposition \ref{pdinbns} extends to $X_Q$ by the K\"unneth formula,
which we express in terms of restriction along $i_Q$ in our
applications below; it also extends to general doubly complex oriented
cohomology theories in the obvious fashion. It inspires many
interesting cobordism calculations, of which we offer a single
example.
\begin{prop}\label{singam}
The map $q_h\colon B_n\rightarrow\CP^{n-h+1}$ represents either of the
expressions  
\[
\sum_{m=0}^{n+1-h}g_{n-m}\beta_{m,\ell}
\sptsp{or}
\sum_{j\geq m=0}^{n+1-h}g_{n-j}(g)^m_{j-m}\beta_{m,r}
\]
in $\varOmega^{DU}_{2n}(\CP^{n-h+1})$, for each $1\leq h\leq n$.
\end{prop}
\begin{proof}
The coefficient of $\beta_{m,\ell}$ in the first
expression is given by $\langle x_{h,\ell}^m,x_{[n]}\rangle$; by
Proposition \ref{pdinbns}, this is $g_{n-m}$ when $1\leq m\leq n-h+1$,
and zero otherwise, as required. To convert the result into the second
expression, we dualize the expansion \eqref{xrxs}. 
\end{proof}

\section{Applications}\label{ap}

In our final section, we apply the duality calculations to realize the
left and right actions of the Landweber--Novikov algebra on its dual;
some preliminary combinatorics is helpful.

Fixing the subset $Q=\cup_{j=1}^sI(j)$ of $[n]$, we consider the
additive semigroup $H(Q)$ of nonnegative integer sequences $h$ of the
form $(h_1,\dots,h_n)$, where $h_i=0$ for all $i\notin Q$; for any such
$h$, we set $|h|=2\sum_i h_i$. Whenever $h$ satisfies
$\sum_{i=l}^{b(j)}h_i\leq b(j)-l+1$ for all $a(j)\leq l\leq b(j)$, we
define the subset $hQ\subseteq Q$ by
\[
\{m:\text{$\sum_{i=l}^mh_i<m-l+1$ 
for all $a(j)\leq l\leq m\leq b(j)$}\};
\]
otherwise, we set $hQ=Q$. It follows that $hQ=Q\cap h[n]$ for all $h$
in $H(Q)$, and we introduce the subset $S(h)\subseteq[s]$ of indices
$j$ for which $I(j)\cap hQ\neq\emptyset$. We also identify the
subsemigroup $K(Q)\subseteq H(Q)$ of sequences $k$ for which $k_i$ is
nonzero only if $i=a(j)$ for some $1\leq j\leq s$.

For each $h$ in $H(Q)$ and $k$ in $K(Q)$, our applications require us
to invest the manifold $X_{Q,(h+k)[n]}$ of Lemma \ref{intersecs} with 
a non-basic double--$U$ structure. In terms of the decomposition
$\dtimes_{S(h+k)}Y_{I(j)\cap(h+k)[n]}$, this is given by 
\begin{equation}\label{newdoubu}
\Big(\dtimes_{S(h+k)}(\gamma-k_{a(j)}\gamma_1),\; 
\dtimes_{S(h+k)}(k_{a(j)}+1)\gamma_1\Big),
\end{equation}
and we denote the resulting double--$U$ manifold by $X^k_{Q,(h+k)[n]}$.
For example, when $h$ is $0$ and $k$ has a single nonzero element
$k_{a(j)}=m$ for some $1\leq j\leq s$ and $m\leq b(j)-a(j)$, then
$X^k_{Q,(h+k)[n]}$ reduces to the manifold
$X_{Q\setminus[a(j),a(j)+m-1]}$ with double $U$--structure 
\begin{equation}\label{xqmj}
\big(
\gamma^{\times j-1}\times(\gamma-m\gamma_1)\times\gamma^{\times s-j},
\;\gamma_1^{\times j-1}\times(m+1)\gamma_1\times\gamma_1^{\times s-j}
\big).
\end{equation}
This case is important enough to motivate the notation
$X^{m:j}_P$ (omitting the $:j$ if $s=1$) for any $X_P$
whose basic double $U$--structure is similarly amended on its $j$th
factor $Y_{I(j)}$; in particular, \eqref{xqmj} describes
$X^{m:j}_{Q\setminus[a(j),a(j)+m-1]}$. 

We may now apply Proposition \ref{pdinbns} to compute the effect of
the left and right actions of $S^*$ on $S_*$ under the canonical
isomorphism. To ease computations with the left action we consider 
the monomial basis of {\it tangential\/} Landweber--Novikov operations
$\bar{s}_\psi$ for $A_*^U$; under the universal Thom
isomorphism, these correspond to the Chern classes $\perp^*\!c_\psi$
induced by the involution $\perp$ of complementation on $\BU$. There are
therefore expressions 
\begin{equation}\label{barexp}
\bar{s}_\psi=\sum_\omega\lambda_{\psi,\omega}s_\omega,
\end{equation}
where the $\lambda_{\psi,\omega}$ are integers and the summation
ranges over sequences $\omega$ for which $|\omega|=|\psi|$ and
$\sum\omega_i\geq\sum\psi_i$. For each $Q\subseteq[n]$, it is also
helpful to partition $K(Q)$ and $H(Q)$ into compatible blocks
$K(Q,\psi)$ and $H(Q,\psi)$ for every indexing sequence $\psi$; each
block consists of those sequences $k$ or $h$ which have $\psi_i$
entries $i$ for each $i\geq 1$, and all other entries zero. Thus, for
example, $|h|=|\psi|$ for all $h$ in $H(Q,\psi)$. Any such block will
be empty whenever $\psi$ is incompatible with $Q$ in the appropriate
sense. 
 
\begin{thm}\label{specss}
Up to double $U$--cobordism, the actions of $S^*_\ell$ and $S^*_r$ on
additive generators of $G_*$ are induced by
\[
\bar{s}_{\psi,\ell}(X_Q)=
\sum_{H(Q,\psi)}X_{Q,h[n]}
\spandsp
s_{\omega,r}(X_Q)=\sum_{K(Q,\omega)}X^k_{Q,k[n]}
\]
respectively.
\end{thm}
\begin{proof}
We combine \eqref{onptdoub} with Proposition \ref{pdinbns}, recalling
that $c_\theta$ is evaluated on any sum of line bundles
$\oplus_{i=1}^r\lambda_i$ by forming the symmetric sum of all
monomials $c_1(\lambda_1)^{i_1}\dots c_r(\lambda_r)^{i_r}$, where
$\theta_i$ of the exponents take the value $i$ for each $1\leq i\leq
r$. We note that the product structure in $\odus(B_n)$ allows us to
replace any $x_i^m$ (either left or right) by $x^{[i,i+m-1]}$ when
$[i,i+m-1]\subseteq Q$, and zero otherwise; indeed, the definitions of
$H(Q)$ and $K(Q)$ are tailored exactly to these relations. For
$s_{\psi,\ell}(X_Q)$ we set $k=0$, and observe that
$\bar{c}_{\psi,\ell}(\nu_\ell)=i_Q^*c_{\psi,\ell}(\oplus_Q\gamma_i)$.
For $s_{\omega,r}(X_Q)$ we set $h=0$, and observe in turn that
$c_{\omega,r}(\nu_r)=
i_Q^*c_{\omega,r}(\gamma_{a(1)}\oplus\dots\oplus\gamma_{a(s)})$. The
computations are then straightforward, although the bookkeeping
demands caution. 
\end{proof}

Recalling \eqref{newdoubu}, we may combine the left and right actions
by 
\[
\bar{s}_{\psi,\ell}\otimes s_{\omega,r}(X_Q)=
\sum_{H(Q,\psi),\;K(Q,\omega)}X^k_{Q,(h+k)[n]},
\]
from which the diagonal action follows immediately. If we prefer to
express the action of $S^*_\ell$ in terms of the standard basis
$s_\omega$, we need only incorporate the integral relations
\eqref{barexp}. 

Readers may observe that our expression in \S\ref{boflma} for
$\nu_\ell$ as the sum of line bundles $\bigoplus_{i=1}^ni\rho_i$
appears to circumvent the need to introduce the tangential operations
$\bar{s}_\psi$. However, it contains $n(n+1)/2$ summands rather
than $n$, and their Chern classes $y_i$ are algebraically more
complicated than the $x_i$ used above, by virtue of \eqref{xitoy}. 
These two factors conspire to make the alternative
calculations less palatable, and it is an instructive exercise to
reconcile the two approaches in simple special cases. The apparent
dependence of Theorem \ref{specss} on $n$ is illusory (and solely for
notational convenience), since $k_i$ and $h_i$ are zero whenever $i$
lies in $Q'$.

We may specialize Theorem \ref{specss} to the cases when 
$\psi$ and $\omega$ are of the form $\epsilon(m)$ for some integer
$0\leq m\leq |Q|$, or when $Q=[n]$ (so that we are dealing with
polynomial generators of $G_*$), or both. We obtain  
\begin{multline}\label{vspec1}
\bar{s}_{\epsilon(m),\ell}(X_Q)=
\sum_j\sum_{i=a(j)}^{b(j)-m+1}X_{Q\setminus I(j)}\times 
Y_{I(j)\setminus[i,i+m-1]}\\
\sands
s_{\epsilon(m),r}(X_Q)=
\sum_jX^{m:j}_{Q\setminus[a(j),a(j)+m-1]},
\end{multline}
where the summations range over all $j$ with $b(j)-a(j)\geq m-1$, and
\begin{multline}\label{vspec2}
\bar{s}_{\psi,\ell}(X_{[n]})= 
\sum_{H([n],\psi)}Y_{h[n]}\\
\sands
s_{\omega,r}(X_{[n]})=
\begin{cases}
X^m_{[m+1,n]}&\text{when $\omega=\epsilon(m)$}\\
0&\text{otherwise}.
\end{cases}
\end{multline}
These follow from \eqref{newdoubu}, and the facts that
$K(Q,\epsilon(m))$ consists solely of sequences containing a single
nonzero entry $m$ in some position $a(j)$, and $K([n],\omega)$ is
empty unless $\omega=\epsilon(m)$ for some $0\leq m\leq n$.

We might expect Theorem \ref{specss} to provide geometrical
confirmation that $G_*$ is closed under the action of $S^*_\ell\otimes
S^*_r$ on $\osdu$, as noted in Proposition \ref{lislnrisr}; however,
it remains to show that $X^{k+1}_{kQ}$ lies in $G_*$! Currently, we
have no direct geometrical proof of this fact.

We now turn to the structure maps of $S_*$, continuing to utilize the
canonical isomorphism to identify $G_*$ and $G_*\otimes G_*$ with $S_*$
and $S_*\otimes S_*$ respectively. We express monomial generators of
$G_*\otimes G_*$ as double $U$--cobordism classes of {\it pairs\/} of
basic double $U$--manifolds $(X_Q,X_R)$, where $Q$ and $R$ range over
independently chosen subsets of $[n]$. 
\begin{prop}\label{mainthm}
Up to double $U$--cobordism, the coproduct $\delta$ and the antipode
$\chi$ of the dual of the Landweber--Novikov algebra are induced by 
\[
X_Q\mapsto\sum_{K(Q)}(X^k_{Q,k[n]},\,X_{Q\setminus kQ})
\spandsp 
X_Q\mapsto\chi(X_Q)
\]
respectively.
\end{prop}
\begin{proof}
For $\delta$, we combine the right action of Theorem \ref{specss} with
\eqref{lractalt}, and the observation that $X_{Q\setminus kQ}$ is
isomorphic to $B^\omega$ for each $k$ in $K(Q,\omega)$. For $\chi$, we
refer to Proposition \ref{bnisbn}.
\end{proof}
\begin{cor}\label{realgk}
When equipped with the double $U$--structure 
\[
\Big(\dtimes_{j=1}^s(\gamma-m(j)\gamma_1),
\;\dtimes_{j=1}^s(m(j)+1)\gamma_1\Big),
\]
the manifold $X_Q$ represents $\prod_j(g)_{b(j)-a(j)+1}^{m(j)+1}$ in
$\varOmega_{2|Q|}^{DU}$ for any sequence of natural numbers $m(1)$,
$m(2)$, \dots, $m(s)$. 
\end{cor}
\begin{proof}
If we consider the coproduct for $Q=[n]$ in Proposition \ref{mainthm},
we deduce that $X^m_{[m+1,n]}$ represents $(g)_{n-m}^{m+1}$ by
appeal to \eqref{lnstrucs}. The result for general $X_Q$ follows by
applying this case to each factor $Y_{I(j)}$.
\end{proof}
Corollary \ref{realgk} is particularly fascinating because it
describes how to represent an intricate (but important) polynomial in
the cobordism classes of the basic $B_n$ by perturbing the double
$U$--structure on a single manifold $X_Q$.  

For a final comment on Proposition \ref{mainthm}, we note that the
elements of $\osdu\otimes_{\varOmega_*^U}\osdu$ may be represented by
{\it threefold} $U$--manifolds. Under the canonical isomorphism, the
coproduct on the Hopf algebroid $A_*^{DU}$ is then induced by mapping
the double $U$--cobordism class of each $(M;\nu_\ell,\nu_r)$ to the
threefold cobordism class of $(M;\nu_\ell,0,\nu_r)$, and the diagonal
on $G_*$ follows by restriction. Theories of multi $U$--cobordism are
remarkably rich, and have applications to the study of iterated
doubles and Adams--Novikov resolutions; we reserve our development of
these ideas for the future.

\end{document}